\documentclass[12pt]{amsart}
\usepackage[top=30truemm,bottom=30truemm,left=25truemm,right=25truemm]{geometry}
\usepackage{txfonts}
\usepackage{mathrsfs}

\usepackage{color}
\usepackage{bm}
\usepackage{amsfonts,amssymb}
\usepackage{dsfont}
\usepackage{amscd}
\usepackage{extarrows}
\usepackage{amsmath}
\usepackage{mathrsfs}
\usepackage{enumerate}
\usepackage{amscd}
\usepackage[all]{xy}
\usepackage[pagebackref,colorlinks]{hyperref}

\newtheorem{thm}{Theorem}[section]
\newtheorem{defn}[thm]{Definition}
\newtheorem{lem}[thm]{Lemma}
\newtheorem{cor}[thm]{Corollary}
\newtheorem{prop}[thm]{Proposition}
\newtheorem{rem}[thm]{Remark}

\newcommand{\be}{\begin{equation}}
	\newcommand{\ee}{\end{equation}}
\newcommand{\bea}{\begin{eqnarray}}
	\newcommand{\eea}{\end{eqnarray}}
\newcommand{\ben}{\begin{eqnarray*}}
	\newcommand{\een}{\end{eqnarray*}}
\newcommand{\bt}{\begin{split}}
	\newcommand{\et}{\end{split}}
\newcommand{\bet}{\begin{equation}}

\begin{document}
		\title{Solution to Hessian type equations with prescribed singularity on compact K\"ahler manifolds}
		
		%
		\author[G. Li]{Genglong Lin}

		\thanks{}
		
\begin{abstract}
Let $(X,\omega)$ be a compact K\"ahler manifold of dimension $n$ and fix  an integer $m$ such that $1\leq m\leq n$. We reformulate most relative pluripotential results of Darvas-DiNezza-Lu's survey \cite{DNL23} to the Hessian setting. As an application, we use a slightly different method and give an characterization of finite energy range of the Hessian operator, which cannot be directly reformulated by \cite{DNL23}.
 
 Given a model potential $\phi$, we also study degenerate complex Hessian equations of the form $(\omega+dd^c \varphi)^m\wedge\omega^{n-m}=F(x,\varphi)\omega^n$. Under some natrual conditions on $F$, we prove that the solution of this type equation has the same singularity type as $\phi$.
\end{abstract}
	\maketitle
\tableofcontents
\section{Introduction}
The complex $m$-Hessian equation has been studied intensively in recent years and has more and more importance in complex analysis, complex geometry and other fields. It can be seen as an interpolation between the classical Poisson equation $(m=1)$ and the complex Monge-Amp\`ere equation $(m=n)$. The Monge-Amp\`ere case has a lot of applications in algebraic geometry such as the celebrated Calabi-Yau's theorem \cite{Yau78} and the existence of K\"ahler-Einstein metric on compact K\"ahler manifolds, specially, the famous YTD conjecture. There are various relaxation of conditions of the equation, see \cite{DNL23}\cite{BEGZ10}\cite{BBGZ13}\cite{EGZ09}, to only cite a few. Some generalization also lead to some conjectures such as Tosatti-Weinkove conjecture and Demailly-Paun's conjecture (see \cite{Ngu16}). For the Hessian case, it appears in the study of the Fu-Yau equation related to the Strominger system \cite{PPZ17}\cite{PPZ18}\cite{PPZ19}, which is motivated by the study of the Calabi problem for HKT-manifolds \cite{AV10}.

The program of solving the non-degenerate complex Hessian equation on compact K\"ahler manifold, to the author's knowledge, started from \cite{Hou09},\cite{Jbi10} and \cite{Kok10}. Their work requires some conditions on the underlined manifold. After \cite{HMW10} provided a general a priori $\mathcal{C}^2$ estimate, \cite{DK17} solved the equation at full generality. They also have given a very powerful $\mathcal{C}^0$-estimate in \cite{DK14}, which allows one to find a continuous weak solution of the degenerate complex Hessian equation when the RHS of the equation has $L^p(x),p>n/m$ density with repsect to the volume form. Now let $(X,\omega)$ be a compact K\"ahler manifold of complex dimension $n$ and fix $1\leq m\leq n$. In this paper we study degenerate complex Hessian equations of the form:
\begin{equation}\label{eqq}
(\omega+dd^c \varphi)^m\wedge\omega^{n-m}=F(x,\varphi)\omega^n
\end{equation}
where $F:X\times\mathbb{R}\to\mathbb{R}^{+}$ satisfies some natrual integrability conditions.
 
Two special cases  are $F(x,t)=f(x)$ and $F(x,t)=f(x)e^t$. When $F(x,t)=f(x)$ and $m=n$, it becomes the usual Hessian equation which corresponds to the usual Monge-Amp\`ere equation. When $F(x,t)=f(x)e^t$ and $m=n$, it becomes the usual Hessian type equation which corresponds to the usual K\"ahler-Einstein equation of $\lambda>0$.
 
To deal with degenerate complex Hessian equaiton, many scholars develop potential theory to explore more and more properties of Hessian operator and the corresponding convergence theorem, see \cite{Blo05}\cite{Lu13}\cite{LN15}\cite{LN22} and references therein. \cite{Blo05} aimed at basic, mostly local properties of
$m$-subharmonic functions and the Hessian operator $H_m$, whereas \cite{Lu13} first defined the class of $(\omega,m)$-subharmonic functions and the Hessian operator for bounded $(\omega,m)$-subharmonic functions which is delicate because of the lack of global regularization process. \cite{Pli13}\cite{LN15} used different methods to solve this global regularization problem for $\omega$-$m$-subharmonic functions. By using variational approach due to Berman-Boucksom-Guedj-Zeriahi\cite{BBGZ13}, \cite{LN15} also solved degenerate complex Hessian equations whose RHS is an arbitrary probability measure (non-$m$-polar measure) which do not charge $m$-polar sets. They also solved equations (\ref{eqq}) when $F(x,t)=f(x)e^t$ and found the solution belonging to the full Hessian mass class $\mathcal{E}(X,\omega,m)$. For the general  case of $F(x,t)$, \cite{Lu13} solved Hessian type equation (\ref{eqq}) whose solution is continuous on $X$. The proof essentially used techniques from \cite{EGZ09}, more precisely, Ko\l odziej's $L^{\infty}$ estimate \cite{Kol98}.

Recently, \cite{DNL23}\cite{DNL21}\cite{DNL18}estabilished relative pluripotential theory and studied Monge-Amp\`ere equations that have a prescribed singularity profile. This type covers the usual Monge-Amp\`ere studied by \cite{BEGZ10}\cite{BBGZ13}. It also has natruality at finding certain quasi-psh functions that have prescribed singularity, for example analytic singularity around given submanifolds.  The results was promoted to the Hessian equation by \cite{LN22} by using the supersolution method of \cite{GLZ19} as in \cite{DNL21}. They considered the case that the RHS of (\ref{eqq}) was either $F(x,t)=f(x)$ (more generally a non-$m$ polar measure) or $F(x,t)=f(x)e^t$. \cite{AAG22} also considered such Hessian type equation in a general situation that the measure does not charge $m$-polar subset and for all $t\in\mathbb{R},x\mapsto F(x,t)\in L^1(\mu)$ and prove the existence of solution. One may ask about the regularity of the solution of this type equation. We here state:

\noindent {\bf Theorem A:}
Assume that $\phi$ is a model potential i.e. $\phi=P[\phi]$ (see definition (\ref{definition of rooftop}) of rooftop envelope of $\omega$-$m$-sh functions) and $\int_X H_m(\phi)>0$. Let $F:X\times\mathbf{R}\rightarrow[0,+\infty)$ be a function satisfying the following conditions:
\begin{itemize}
\item[(i)] for almost every $x\in X, t\to F(x,t)$ is non-decreasing and continuous,
\item[(ii)] for any fixed $t\in\mathbf{R}$ there exists $p>n/m$ such that the function $x\to F(x,t)$ belongs to $L^p(X)$,
\item[(iii)] there exists $t_0\in\mathbf{R}$ such that $\int_X F(\cdot,t_0)=\int_X \omega^n$.
\end{itemize}
Then the soution of above equation $\varphi\in\mathcal{E}_{\phi}(X,\omega,m)$ (relative full Hessian mass class) has the same singularity type as $\phi$.
\begin{rem}
By the classical philosophy of solving Monge-Amp\`ere or Hessian equation, it seems that (relative) boundedness of solution requires RHS of the equation to have $L^p,p>1$ density with respect to the volume form. With such consideration, we keep the original situation of \cite{Lu13} and give an alternative proof of existence theorem which is a special case of \cite{AAG22}. 
\end{rem}

In the latest Darvas-DiNezza-Lu's survey \cite{DNL23} about relative pluripotential theory on compact K\"ahler manifold, they gave an general answer of Guedj-Zeriahi's question \cite{GZ07} about the finite energy range of Monge-Amp\`ere operator. Their proof used Skoda's integrability theorem, which has no parallel version of quasi-$m$-sh function. Motivated by \cite{LN15}, we use a slightly different approach to skip this obstruction. Namely, we also give an characterization of finite energy range of Hessian operator on compact K\"ahler manifolds:

\noindent {\bf Theorem B:}
Assume that $\phi$ is a model potential and $\mu$ is a Radon measure such that $\mu(\{\phi=-\infty\})=0,\int_X H_m(\phi)=\mu(X)>0$. Then the following statements are equivalent:
\begin{itemize}
\item[(i)] There exists a constant $C>0$ such that for all $u\in\mathcal{E}_{\chi}(X,\omega,m,\phi)$ with $\sup_X u=0$, we have
$$\int_{X}\chi(\phi-u)d\mu\leq CE_{\chi}(u,\phi)^{M/(M+1)}+C.$$
\item[(ii)]$\chi(|\phi-u|)\in L^{1}(\mu),\forall u\in\mathcal{E}_{\chi}(X,\omega,m,\phi)$.
\item[(iii)]$\mu=H_m(\varphi)$ for some $\varphi\in\mathcal{E}_{\chi}(X,\omega,m,\phi)$ satisfying $\sup_X \varphi=0$,

\end{itemize}
where $\chi:[0,\infty)\to[0,+\infty)$ is a continous increasing function such that
$\forall t\geq 0,\forall \lambda\geq 1, \chi(\lambda t)\leq\lambda^M\chi(t),1\leq M<\frac{n}{n-m}$.

\section{Preliminaries}
Let $(X,\omega)$ be a compact K\"ahler manifold of dimension n, and fix an integer $m$ such that $1\leq m\leq n$.
Fix $\Omega$ an open subset of $\mathbf{C}^n$ and $\beta:=dd^c \rho$ a K\"ahler form in $\Omega$ with smooth bounded potential.
\begin{defn}
A function $u\in C^2(\Omega,\mathbf{R})$ is called $m$-subharmonic ($m$-sh for short) with respect to $\beta$ if the following inequalities hold in $\Omega$:
$$(dd^cu)^k\wedge\beta^{n-k}\ge0,\forall k\in\{1,...,m\}.$$
\end{defn}
\begin{defn}
A function $u\in L^1(\Omega,\mathbf{R})$ is called $m$-subharmonic with respect to $\beta$ if
\begin{itemize}
\item[(1)] $u$ is strictly upper semicontinuous in $\Omega$, namely, $\limsup_{z\to z^0}u(z)=u(z_0)$ for any $z_0\in\Omega$,
\item[(2)] $dd^cu\wedge dd^cu_2\wedge...\wedge dd^cu_m\wedge\beta^{n-m}\geq0$, for all $u_{2},...,u_{m}\in C^{2}(\Omega)$, $m$-sh with respect to $\beta$.
\end{itemize}
\end{defn}
As observed by \cite{DK14}, the two definitions of $m$-sh functions above coincide for smooth functions.
\begin{defn}
A function $u\in L^1(X,\omega^n)$ is called $\omega$-$m$-subharmonic ($\omega$-$m$-sh for short)if, locally in $\Omega\subset X$ where $\omega=dd^c \rho, u+\rho$ is $m$-subharmonic with respect to $\omega$.
The set of all $\omega$-$m$-sh functions on $X$ is denoted by $\mbox{SH}_m(X,\omega)$.
\end{defn}
One can locally regularize $m$-sh functions by using smooth convolution kernel. As for the global regularization for $\omega$-$m$-sh functions on  compact K\"ahler manifold $X$, this is more complicated but valid by \cite{LN15}.

Given $u,v\in\mbox{SH}_m(X,\omega)$, similarly to quasi-psh functions, we say that $u$ is less singular than $u$ if there exists a constant $C$ such that $v\leq u+C$. We say that $u$ has the same singularity type as $v$ if there exists a constant $C$ such that $u-C\leq v\leq u+C$.
In the flat case, Blocki proved in \cite{Blo05} that $m$-sh functions are in $L^p$ for any $p<n/(n-m)$, and conjectured that it holds for $p< nm/(n-m)$. Using the $L^{\infty}$ estimate due to Dinew and Kolodziej, one can prove the same integrability property for $\omega$-$m$-sh functions, see \cite{LN15}.
\begin{defn}
Following Bedford-Taylor (see also \cite{LN15}), we can define complex hessian operator for any $\omega$-$m$-sh functions $u_1,u_2,...,u_m$ (not necessarily bounded):
$$H_m(u_1,u_2,...,u_m):=\lim_{s\to {+\infty}}\mathbf{1}_{U^s}(u_1^s,...,u_m^s)$$
where $U^{s}:=\cap_{p=1}^{m}\{u_{p}>-s\}$ and $u^{s}:=\max(u,-s)$.

\end{defn}
When $u_1=...=u_m=u$ we simply denote the Hessian measure $H_m(u,u,...,u)$ by $H_m(u)$. By stokes' theorem,
$$0\leq\int_{X}H_{m}(u)\leq1.$$
A Borel set $E$ is called $m$-polar (with respect to $\omega$) if there exists $u\in\mbox{SH}(X,\omega)$ such that $E\subset\{u=-\infty\}$. Since $H_k(u)$ is the strong limit of $\mathbf{1}_{\{u>-j\}}H_k(u_j)$, the positive measure $H_k(u)$ does not charge $m$-polar sets (see \cite[Lemma 2.4]{LN22}).

For a Borel set $E\subset X$ we define 
$$\operatorname{Cap}_{m}(E):=\sup\left\{\int_{E}H_{m}(u)\mid u\in\operatorname{SH}_{m}(X,\omega),-1\leq u\leq0\right\}.$$
A sequence of functions $u_j$ converges in capacity to $u$ if for all $\epsilon>0$,
$$\lim\limits_{j\to+\infty}\mathrm{Cap}_{m}(|u_{j}-u|>\varepsilon)=0.$$
\begin{defn}
A Borel set $E\subset X$ is called quasi-open (quasi-closed) if for each $\epsilon>0$, there exists an open (closed) set $U$ such that
$$\operatorname{Cap}_m((E\setminus U)\cup(U\setminus E))<\varepsilon.$$
\end{defn}
We recall a classical convergence theorem from Bedford-Taylor theory. By a similar proof to \cite[Theorem 4.26]{GZ17}, we have:
\begin{prop}
Let $U\subset\mathbf{C}^n$ be an open set. Suppose $\{f_j\}_j$ are uniformly bounded quasi-continuous functions which converge in capacity to another quasi-continuous function $f$ on $U$. Let $\{u_1^j\}_j,...,\{u_m^j\}_j$ be uniformly bounded $m$-sh functions on $\Omega$, converging in $m$-capacity to $u_1,...,u_m$ respectively. Then we have the following weak convergence of measures
$$f_ji\partial\bar{\partial}u_1^j\wedge i\partial\bar{\partial}u_2^j\wedge\ldots\wedge i\partial\bar{\partial}u_m^j\wedge\beta^{n-m}\rightarrow fi\partial\bar{\partial}u_1\wedge i\partial\bar{\partial}u_2\wedge\ldots\wedge i\partial\bar{\partial}u_m\wedge\beta^{n-m}$$
where $\beta$ is the standard K\"ahler form of $\mathbf{C}^n$
\end{prop}
Applying above fact, one can easily deduce that
\begin{lem}
Assume $u_j$ is a sequence of uniformly bounded $\omega$-$m$-sh functions in $U\subset X$ converging in $m$-capacity to $\omega$-$m$-sh function $u$. Suppose $\{f_j\}$ are uniformly bounded non-negative quasi-continuous functions which converge in $m$-capacity to another quasi-continuous function $f\geq0$ on $U$.

If $E\subset U$ is a quasi-open set then
$$\liminf_{j\to\infty}\int_{E}f_j(\omega+dd^cu_j)^n\geq\int_{E}f(\omega+dd^cu)^n.$$
If $V\subset U$ is a quasi-closed set then
$$\operatorname*{lim}_{j\to\infty}\int_{V}f_{j}(\omega+dd^{c}u_{j})^{n}\leq\int_{V}f(\omega+dd^{c}u)^{n}.$$
\end{lem}

The following lower-semicontinuity property of Hessian product will be key in the sequel.
\begin{thm}\cite[Theorem 3.3]{LN22}\cite[Theorem 2.6]{DNL23}\label{convergence thm}
Assume that $u_1^j,...,u_m^j$ are sequences of $\omega$-$m$-sh functions converging in $m$-capacity to $\omega$-$m$-sh functions $u_1,...,u_m$. Let $\chi_j$ be a sequence of positive uniformly bounded quasi-continuous functions which converges in $m$-capacity to $\chi$. Then
$$\liminf_{j\to+\infty}\int_{X}\chi_{j}H_{m}(u_{1}^{j},...,u_{m}^{j})\geq\int_{X}\chi H_{m}(u_{1},...,u_{m}).$$
If additionally,
$$\int_X H_m(u_{1},...,u_{m})\geq\limsup_{j\to{+\infty}}\int_X H_m(u_{1}^{j},...,u_{m}^{j}),$$
then
$$\chi_j H_m(u_1^j,...,u_m^j)\to\chi H_m(u_1,...,u_m)$$ in the sense of measures.
\end{thm}
\section{Relative m-potential theory}
\subsection{Quasi-$m$-sh envelopes}
If $f$ is a function on $X$, we define the envelope of $f$ in the class $\mbox{SH}_m(X,\omega)$ by
$$P(f):=(\sup\{u\in\mbox{SH}(X,\omega):u\leq f\})^{*},$$
with the convention that $\sup\emptyset=-\infty$. Notice that $P(f)\in\mbox{SH}(X,\omega)$ if and only if there exists some $u\in\mbox{SH}_m(X,\omega)$ lying below $f$. Note also that $P(f+C)=P(f)+C$ for any constant $C$. If in particular $f=\min(\psi,\phi)$, we denote the envelope by $P(\psi,\phi):=P(\min(\psi,\phi))$. Oberve that $P(\psi,\phi)=P(P(\psi),P(\phi))$, so w.l.o.g. we can assume that $\psi,\phi$ are two $\omega$-$m$-sh functions. We recall and show some important properties about envelopes:
\begin{thm}\label{concentr thm1}
Assume that $f$ is quasi-continuous and usc on $X$. Then
$$\int_{\{P(f)<f\}}H_m(P(f))=0.$$
\end{thm}
\begin{proof}
Without loss of generality we can assume that $f\leq0$. By \cite[Corollary 3.4]{LN15},if $f\in C(X)$, then the hessian measure $H_m(P(h))$ vanishes on $\{P(h)<h\}$. To treat the general case, we can adapt the idea of \cite[Proposition 2.16]{DNL18} since we can approximate $f$ from above by a sequence of smooth function $(f_j)$. We omit the proof here.
\end{proof}

\cite{LN22} also proves a mass concentration theorem when $f$ is of a special type:
\begin{thm}\cite[Proposition 3.10]{LN22}\label{concentr thm2}
Assume that $f=a\varphi-b\psi$, where $\varphi,\psi\in\mbox{SH}_m(X,\omega)$, and $a,b$ are positive constants. If $P(f)\not\equiv-\infty $ then
$$\int_{\{P(f)<f\}}H_{m}(P(f))=0.$$
Here the function $f=a\varphi-b\psi$ is well-defined in the complement of a m-polar set.
\end{thm}
\begin{rem}
By the proof of \cite[Theorem 2.7]{DNL23}, we can prove that if $f$ is quasi-continuous on $X$ and $P(f)\in\mbox{SH}_m(X,\omega)$, $H_m(P(f))$ does not charge the subset $\{P(f)=f\}$. This generalizes the above result.
\end{rem}

\begin{lem}\cite[Proposition 2.10, Lemma 2.9]{LN22}\cite[Lemma 2.9]{DNL23}\label{compare measure}
Let $\varphi,\psi\in\mbox{SH}(X,\omega)$. Then
$$H_m(\max(\varphi,\psi))\geq\mathbf{1}_{\{\psi\leq\varphi\}}H_m(\varphi)+\mathbf{1}_{\{\varphi<\psi\}}H_m(\psi).$$
In particular, if $\varphi\leq\psi$ then $\mathbf{1}_{\{\varphi=\psi\}}H_m(\varphi)\leq\mathbf{1}_{\{\varphi=\psi\}}H_m(\psi)$.
\end{lem}
In order to investigate the relative $m$-potential theory, starting form the "rooftop envelope" we introduce the 
\begin{equation}\label{definition of rooftop}
P[\psi](\phi):=\left(\lim_{C\to\infty}P_\omega(\psi+C,\phi)\right)^*.
\end{equation}
It is easy to see that $P[\psi](\phi)$ only depends on the singularity type of $\psi$. When $\phi=0$, we simply
write $P[\psi]=P[\psi](0)$, and we refer to this potential as the envelope of the singularity type $[\psi]$.
\begin{lem}\label{100}
The operator $P[\cdot](\phi)$ is concave: if $u,v,\phi\in\mbox{SH}(X,\omega)$ and $t\in(0,1)$ then
$$P[tu+(1-t)v](\phi)\geq tP[u](\phi)+(1-t)P[v](\phi).$$
\end{lem}

\begin{prop}\cite[Corollary 3.11]{LN22}\cite[Theorem 3.6]{DNL23}\label{104}
Let $u,v\in\mbox{SH}_m(X,\omega)$ be such that $P(u,v)\in\mbox{SH}_m(X,\omega)$. Then
$$H_{m}(P(u,v))\leq\mathbf{1}_{\{P(u,v)\boldsymbol{=}u\}}H_{m}(u)+\mathbf{1}_{\{P(u,v)\boldsymbol{=}v\}}H_{m}(v).$$ In particular, $H_{m}(P[u])\leq\mathbf{1}_{\{P[u]=0\}}\omega^{n}.$ $H_m(P[\psi](\varphi))\leq\mathbf{1}_{\{P[\psi](\varphi)=\varphi\}}H_m(\varphi)$.
\end{prop}

\subsection{Monotonicity of Hessian product masses}
Following the idea of \cite[\S 3]{DNL23}, we give another proof of the monotonicity of Hessian masses. For $m=n$ the Monge-Amp\`ere case, the result is due to \cite{Wit19}. Using the monotonicity of the energy functional, \cite{LN22} prove the result.
\begin{lem}\label{mass eq}
Let $u,v\in\mbox{SH}_m(X,\omega)$. If $u$ and $v$ have the same singularity type, then $\int_X H_m(u)=\int_X H_m(v)$.
\end{lem}
\begin{proof} 
We only sketch the proof. The proof has two steps. Firstly, if there exists a constant $C>0$ such that $u=v$ on the open set $U:=\{\min(u,v)<-C\}$ and $t>C$ fixed, we have $H_m(\max(u,-t))=H_m(\max(v,-t))$ on $U$. Observing that $\{u\leq -t\}=\{v\leq -t\}\subset U$ and $\int_X H_m(\max(u,-t))=\int_X H_m(\max(v,-t))=\mbox{Vol}{(\omega)}$, we have that 
$$\int_{\{u>-t\}}H_m(u)=\int_{\{v>-t\}}H_m(v).$$
Let $t\to\infty$, the first step finishes.

To treat the general case, we may assume that $v\leq u\leq v+B\leq0$ for some positive constant $B$. For each $a\in(0,1)$ we set $v_a:=av,u_a:=\max(u,v_a)$ and $C:=Ba(1-a)^{-1}$. It is easy to check that $u_a=v_a$ on the open set $U_a:=\{\min(u_a,v_a)<-C\}$. Applying the first step we get $\int_X H_m(u_a)=\int_X H_m(v_a)$ as $a\to1$. Expanding $H_m(v_a)$ by linearity of Hessian operator, we get $\int_X H_m(v_a)\to\int_X H_m(v)$. But Proposition \ref{convergence thm} imples 
$\liminf_{a\to 1^{-}}\int_X H_m(u_a)\geq\int_X H_m(u).$ As a result,  $\int_X H_m(u)\leq\int_X H_m(v)$. Exchanging $u$ and $v$, the result holds.
\end{proof}
Applying the above lemma, we also have 
\begin{thm}\cite[Proposition 3.2]{DNL23}\cite[Theorem 3.7]{LN22}\label{mass ineq}
Let $u_1,...,u_m,v_1,...,v_m\in\mbox{SH}_m(X,\omega)$ and assume that $u_j$ is more singular than $v_j$ for all $j$. Then
$$\int_XH_m(u_1,...,u_m)\leq\int_XH_m(v_1,...,v_m).$$
\end{thm}
\begin{rem}\label{P-env mass}
If $u_j^k\nearrow u_j$ a.e. as $k\to\infty$,  then $u_j^k\to u_j$ in capacity by a similar argument of \cite[Proposition 4.25]{GZ17}, and by Proposition \ref{convergence thm} we have $\int_X (\omega+dd^c u_1)\wedge...\wedge(\omega+dd^c u_m)\wedge\omega^{n-m}\geq\limsup_k\int_X (\omega+dd^c u_1^k)\wedge...\wedge(\omega+dd^c u_m^k)\wedge\omega^{n-m}$.
Since for all $u\in\mbox{SH}_m(X,\omega)$ we have $P(u,C)\nearrow P[u]$ as $C\to\infty$, we get
$$\int_X H_m(u)=\int_X H_m(P[u]).$$
Similarly, for $u_j\in\mbox{SH}(X,\omega^j)$ we have $\int_X(\omega+dd^c u_1)\wedge...\wedge(\omega+dd^c u_m)\wedge\omega^{n-m}=\int_X(\omega+dd^c P[u_1])\wedge...\wedge(\omega+dd^c P[u_m])\wedge\omega^{n-m}.$

\end{rem}
\subsection{Model potentials and relative full mass classes}
\begin{defn}
Given $\phi\in\mbox{SH}_m(X,\omega)$, the relative full mass class $\mathcal{E}_{\phi}:=\mathcal{E}_{\phi}(X,\omega,m)$ is the set of all $\omega$-$m$-sh functions $u$ such that $u$ is more singular than $\phi$ and $\int_X H_m(u)=\int_X H_m(\phi)$.
\end{defn}
\begin{defn}
 A model potential is a $\omega$-$m$-sh function $\phi$ such that $P[\phi]=\phi$ and $\int_X H_m(\phi)>0$.
\end{defn}

We can show the maximality of model potentials by a similar argument of \cite[Theorem 3.14]{DNL23}, precisely, the following set of potentials has a maximal element:
$$F_\phi:=\begin{cases}v\in\operatorname{SH}_m(X,\omega):\phi\leq v\leq0~\text{and}~\int_XH_m(v)=\int_XH_m(\phi)\\\end{cases}\Biggr\}.$$
As shown below, if $\int_X H_m(\phi)>0$, this indeed the case, moreover this maximal potential is equal to $P[\phi]$.
\begin{thm}\label{103}
Assume that $\phi\in\mbox{SH}_m(X,\omega)$ and $\int_X H_m(\phi)>0$. Then
$$P[\phi]=\sup_{v\in F_{\phi}}v.$$
In particular, $P[\phi]=P[P[\phi]]$.
\end{thm}

Using this result, we can characterize $\mathcal{E}_{\phi}$ in the following way which generalizes \cite[Theorem 3.15]{DNL23}:

\begin{thm}\label{106}
Suppose $\phi\in\mbox{SH}_m(X,\omega)$ with $\int_X H_m(\phi)>0$ and $\phi\leq0$. The following are equivalent:
\begin{itemize}
\item[(i)] $u\in\mathcal{E}_{\phi}.$
\item[(ii)] $\phi$ is less singular than $u$, and $P[u](\phi)=\phi$.
\item[(iii)] $\phi$ is less singular than $u$, and $P[u]=P[\phi]$.
\end{itemize}
\end{thm}

\begin{cor}
Suppose $\phi\in\mbox{SH}_m(X,\omega)$ such that $\int_X H_m(\phi)>0$. Then $\mathcal{E}_{\phi}(X,\omega,m)$ is convex. Moreover, given $\psi_1,...,\psi_n\in\mathcal{E}_{\phi}(X,\omega,m)$ we have
$$\int_X (\omega+dd^c \psi_1)^{s_1}\wedge...\wedge(\omega+dd^c \psi_n)^{s_m}\wedge\omega^{n-m}=\int_X H_m(\phi).$$
where $s_j\geq0$ are integers such that $\sum_{j=1}^{n}s_j=n$.
\end{cor}
\begin{proof}
The proof is similar to \cite[Corollary 3.16]{DNL23}.
\end{proof}
By a similar proof of \cite[Proposition 3.22]{DNL23} the comparison principle for functions of relative full  Hessian mass also holds:
\begin{prop}\label{general compare principle}
Suppose $\psi_k\in\mbox{SH}_m(X,\omega),k=1,...,j\leq m$ and $\phi\in\mbox{SH}_m(X,\omega)$ is a model potential. If $u,v\in\mathcal{E}_{\phi}$ then
$$\int_{\{u<v\}}\omega_v^{m-j}\wedge\omega_{\psi_1}\wedge...\wedge\omega_{\psi_j}\wedge\omega^{n-m}\leq\int_{\{u<v\}}\omega_u^{m-j}\wedge\omega_{\psi_1}\wedge...\wedge\omega_{\psi_j}\wedge\omega^{n-m}.$$
\end{prop}

\begin{cor}\label{compare principle}
Suppose $\phi\in\mbox{SH}_m(X,\omega)$ is a model potential. If $u,v\in\mathcal{E}_{\phi}$ then
$$\int_{\{u<v\}}H_m(v)\leq\int_{\{u<v\}} H_m(u)\quad\mbox{ and}\quad~\int_{\{u\leq v\}}H_m(v)\leq\int_{\{u\leq v\}}H_m(u).$$
\end{cor}

\subsection{Integration by parts}
The integration by parts formula of Monge-Amp\`ere type was recently obtained \cite{Xia19} using Witt-Nystr\"om’s construction. Adapting ideas of \cite{Lu21}\cite[\S 4]{DNL23} we can prove the integration by parts formula of Hessian type, generalizing their results.
\begin{lem}\label{CLN}
Let $u,v,\psi\in\mbox{SH}_m(X,\omega)$. Assume that $v\leq u\leq v+B$ for some positive constant $B$. Then
$$\int_X \psi H_m(u)\geq\int_X \psi H_m(v)-mB\int_X \omega^n.$$
\end{lem}

\begin{defn}
Given $\phi\in\mbox{SH}_m(X,\omega)$ and $E\subset X$ a Borel subset we define
$$
{\rm Cap}_{m,\phi}(E):=\sup\left\{\int_E H_m(\varphi):\phi-1\leq\varphi\leq\phi\;\text{on}\;\text{X}\right\}
$$
\end{defn}
\begin{rem}\label{Cap m characterize m polar}
From \cite[Lemma 4.3]{DNL18} we know $\mbox{Cap}_{n,\phi}$ characterizes pluripolar set (See also \cite[Lemma 4.9]{LN15} for $\phi=0$). But the relative Hessian case still holds in a similar manner.
\end{rem}

\begin{prop}\label{107}
Assume $\phi \in {\rm PSH}(X,\omega)$ is a model potential.  There exists a constant $C>0$ depending on $X,\omega,n$ such that, for all Borel set $E$, we have 
	$$
	{\rm Cap}_{m,\phi}(E) \leq C {\rm Cap}_{m}(E)^{1/m}. 
	$$
\end{prop}
\begin{proof}
Note that by a standard balayage argument, we have $H_m(V_K^{*})$ is concentrated on $K$ since we can locally solve the Dirichlet problem on any small ball (see \cite{Pli13}). Then adapt the proof of \cite[Proposition 4.4]{DNL23}.
\end{proof}
Follow the line of \cite{DNL23} we can also prove the integration by parts formula, extending \cite[Theorem 3.3]{Lu13} which only applies to bounded $\omega$-$m$-sh functions.

\begin{thm}\label{int by parts}
Let $u,v \in L^\infty(X)$ be differences of quasi-$m$-sh functions, and $\phi_j \in \textup{SH}_m(X,\omega)$, $j \in  \{1,...,m-1\}$. Then
$$\int_X u dd^c v  \wedge \omega_{\phi_1} \wedge \ldots \wedge \omega_{\phi_{m-1}}\wedge\omega^{n-m}= \int_X v dd^c u \wedge \omega_{\phi_1} \wedge \ldots \wedge \omega_{\phi_{m-1}}\wedge\omega^{n-m}.$$
\end{thm}

\section{Complex Hessian type equations with prescribed singularity}
Let $\omega$ be a K\"ahler form on $X$ and $\phi\in\mbox{SH}_m(X,\omega)$.  By $\mbox{SH}_m(X,\omega,\phi)$ we denote the set of $\omega$-$m$-sh functions that are more singular than $\phi$. We say that $v\in\mbox{SH}_m(X,\omega)$ has relatiively minimal singularity type if $v$ has the same singularity type as $\phi$.

Our aim is to consider the following equation of Hessian type:
\begin{equation}\label{hessian type eq}
(\omega+dd^c \phi)^m\wedge\omega^{n-m}=F(x,\phi)\omega^n
\end{equation}
where $F:X\times\mathbf{R}\rightarrow[0,+\infty)$ be a function satisfying the following conditions:
\begin{itemize}
\item[(i)] for almost every $x\in X, t\to F(x,t)$ is non-decreasing and continuous,
\item[(ii)] for any fixed $t\in\mathbf{R}$ there exists $p>n/m$ such that the function $x\to F(x,t)$ belongs to $L^p(X)$,
\item[(iii)] there exists $t_0\in\mathbf{R}$ such that $\int_X F(\cdot,t_0)=\int_X \omega^n$.
\end{itemize}

\cite{AAG22} considered a general situation that the measure does not charge $m$-polar subset and for all $t\in\mathbb{R},x\mapsto F(x,t)\in L^1(\mu)$. But by the classical philosophy of solving Monge-Amp\`ere or Hessian equation, it seems that (relative) boundedness of solution requires RHS of the equation to have $L^p,p>1$ density with respect to the volume form. With such consideration, we keep the original situation of \cite{Lu13} and give an alternative proof of special case of \cite{AAG22}. As an addtional product, we obtain the relative boundedness regularity of solution of above Hessian type equation, which is slightly new.

\subsection{The relative finite energy class}

For $u\in\mathcal{E}_{\phi}$ with relatively minimal singularity type, we define the  Hessian energy of $u$ relative to $\phi$ as
$$\mathrm{I}_{\phi} (u) :=\frac{1}{m+1} \sum_{k=0}^m \int_X (u-\phi) \omega_u^k \wedge\omega_{\phi}^{m-k}\wedge\omega^{n-m}$$
where $\omega_u:=\omega+dd^c u$. The Hessian energy operator enjoys basic good properties like the Monge-Amp\`ere case:

\begin{thm}
 Suppose $u,v \in \mathcal E(X,\omega,\phi)$ have relatively minimal singularity type. Then:
\begin{itemize}
\item[(i)] $ \mathrm{I}_{\phi}(u)-\mathrm{I}_{\phi}(v) = \frac{1}{n+1}\sum_{k=0}^m \int_X (u-v) \omega_{u}^k \wedge \omega_{v}^{m-k}\wedge\omega^{n-m}.$ In particular $\mathrm{I}_{\phi}(u)\leq \mathrm{I}_{\phi}(v)$ if $u\leq v$.
\item[(ii)] If $u\leq \phi$ then, $
\int_X (u-\phi) \omega_u^m\wedge\omega^{n-m} \leq I_{\phi}(u) \leq \frac{1}{n+1} \int_X (u-\phi) \omega_{u}^m\wedge\omega^{n-m}. $ 
\item[(iii)] $\mathrm{I}_{\phi}$ is concave along affine curves. Also, the following estimates hold: $$
	\int_X (u-v) \omega_u^m\wedge\omega^{n-m} \leq I_{\phi}(u) -I_{\phi}(v) \leq \int_X (u-v) \omega_v^m\wedge\omega^{n-m}.$$
\end{itemize}
\end{thm}
\begin{proof}
The proof is essentially achieved by using integration by parts (Theorem \ref{int by parts}). We refer the reader to \cite[Theorem 5.3]{DNL23}.
\end{proof}

We  define the Hessian energy for arbitrary $u\in \mathrm{PSH}(X,\omega,\phi)$ like the Monge-Amp\`ere case:
$$\mbox{I}_{\phi}(u) : =\inf \{\mbox{I}_{\phi}(v) | v\in \mathcal{E}_{\phi}, \; v\ \textrm{has relatively minimal singularity type, and } u\leq v\}. $$
Let also $\mathcal{E}_{\phi}^1(X,\omega,m)$ be the set of all $u\in\mbox{SH}_m(X,\omega,\phi)$ such that $I_{\phi}(u)<+\infty$.

Adapting the ideas of Darvas-DiNezza-Lu, we can alternatively prove the following theorem using variational method. For $\phi=0$ the result is due to \cite{LN22}. See also the supersolution method in \cite{LN22}.
\begin{thm}
Assume that $\mu$ is a positive non-$m$-polar measure on $X$. Then there exists a unique (up to a constant) $u\in \mathcal{E}_{\phi}$ such that 
\begin{equation}\label{eq: MA_exp_version}
	H_m(u) =\mu. 
\end{equation}
\end{thm}
\begin{thm}
Assume that $\mu$ is a positive non-$m$-polar measure on $X$. Then there exists a unique (up to a constant) $u\in \mathcal{E}_{\phi}^1$ such that 
\begin{equation}\label{eq: MA_exp_version}
	H_m(u) =e^{u}\mu. 
\end{equation}
\end{thm}

\subsection{Relative boundedness of solution}
Recall that we work with $\phi \in \mbox{SH}_m(X,\omega)$ such that $P[\phi]=\phi$, and $\int_X H_m(\phi) >0$. Let $f \in L^p(\omega^n),p>n/m$ with $f \geq 0$. In the previous section we have shown that the equation
\begin{equation*}
H_m(u) = f \omega^n, \ \ u \in \mathcal {E}_{\phi}
\end{equation*}
has a unique solution. In this section we will show that this solution has the same singularity type as $\phi$.    Our argument follow the one in \cite{DNL18} which builds on fundamental work of Ko\l odziej in the K\"ahler case (see \cite{Kol98}\cite{Kol03}). We do not follow the one in \cite[\S 5.3]{DNL23} because of a lack of integrability of $e^{-u}$ for $u\in\mbox{SH}_m(X,\omega)$.
\begin{lem}\label{V geq Cap}
Let $1<\beta<\frac{n}{n-m}$. There exists a constant $C=C(p,\omega)$ such that for every Borel subset $K$ of $X$, we have
$$V(K)\leq C\cdot{\rm Cap}_{m,\phi}(K)^{\beta},$$
where $V(K):=\int_K \omega^n$.
\end{lem}
\begin{proof}
The proof is inspired by \cite{DK14}\cite[Lemma 6.6]{LN15}.
Fix an open subset $U$ such that $K\subset U$. Solve the complex Monge-Amp\`ere equation to find $u\in\mathcal{E}_{\phi}$ such that $\omega_u^n=f\omega^n,\sup_X u=0$ on $X$ with $f=V(U)^{-1}\chi_U$. From \cite[Theorem 4.32]{DNL18}, for each $r>1$,
$$0\leq\sup_X (\phi-u)\leq C\|f\|_r^{1/n},$$
where the constant $C$ does not depend on $K$. The inequality between mixed complex Monge-Amp\`ere measures (\cite[Proposition 1.11]{BEGZ10}) tells us that
$$\omega_{u}^{m}\wedge\omega^{n-m}\geqslant f^{m/n}\omega^{n}.$$
Without loss of generality, we can assume that $V(U)<1$. Setting $\lambda=\frac{1}{C\|f\|_r^{1/n}}<1$, $\phi-1\leq\lambda u+(1-\lambda)\phi\leq\phi$, hence we get
$$\mbox{Cap}_{m,\phi}(U)\geq\int_U H_m(\lambda u+(1-\lambda\phi))\geq\lambda^m\int_U H_m(u)\geq C^{-m}V(U)^{1-\frac{m}{rn}}.$$
Thus, for every $r>1$ there exists a constant $C$ not depending on $K$ such that $V(K)\leq C{\rm Cap}_{m,\phi}(K)^{\frac{nr}{nr-m}}$. The proof is complete.

\end{proof}

\begin{prop}\label{vol control cap}
Let $f\in L^p(\omega^n),p>n/m$ with $f\geq0$. Then there exists $C>0$ depending only on $\omega,p$ and $\|f\|_{L^p}$ such that
$$\int_E f \omega^n\leq C{\rm Cap}_{m,\phi}(E)^{1+\alpha}$$
for all Borel sets $E\subset X$.
\end{prop}
\begin{proof}
By Holder's inequality and Lemma \ref{V geq Cap}
$$\int_E f\omega^n\leq\|f\|_{L^p}V(U)^{\frac{p-1}{p}}\leq C\cdot {\rm Cap}_{m,\phi}(E)^{1+\alpha}$$
where $\alpha$ can be  arbitrarily taken from $(0,\frac{mp-n}{(n-m)p})$ because $r$ varies from 1 to $+\infty$.

\end{proof}
\begin{lem}\label{1000}
Let $\phi$ be a model potential and $u,v\in\mathcal{E}_{\phi}$ be two negative functions. Then for all $t>0$ and $\delta\in (0,1]$ we have
$$\operatorname{Cap}_{m,\phi}\{u-v<-t-\delta\}\leq\frac{1}{\delta^{m}}\int_{\{u-v<-t+\delta(\phi-v)\}}H_m(u).$$
\end{lem}
\begin{proof}
Let $\psi\in\mbox{SH}_m(X,\omega,\phi)$ be such that $\phi\leq\psi\leq\phi+1$. In particular $\psi\in\mathcal{E}_{\phi}$. We then have
$$\{u<v-t-\delta\}\subset\{u<(1-\delta)v+\delta\psi-t-\delta\}\subset\{u-v<-t+\delta(\phi-v)\}.$$
Since $\delta^mH_m(\psi)\leq H_m(\delta\psi+(1-\delta)v)$, $u,v$ has relative full mass and $\mathcal{E}_{\phi}$ is convex, Corollary \ref{compare principle} yields
\begin{align*}
\delta^{m}\int_{\{u<v-t-\delta\}}H_m(\psi)& \leq\quad\int_{\{u<\delta\psi+(1-\delta)v-t-\delta\}}H_m(\delta\psi+(1-\delta)v)  \\
&\leq\quad\int_{\{u<\delta\psi+(1-\delta)v-t-\delta\}}H_m(u)\leq\int_{\{u<v-t+\delta(\phi-v)\}}H_m(u).
\end{align*}
Since $\psi$ is an arbitrary candidate in the definition of ${\rm Cap}_{m,\phi}$, the proof is complete.
\end{proof}
\begin{thm}\label{relative boundedness}
Suppose $\phi=P[\phi]$ and $\int_X H_m(\phi)>0$. Let also $\psi\in\mathcal{E}_{\phi}$ with $\sup_X \psi=0$. If $H_m(\psi)=f\omega^n$ for some $f\in L^p(\omega^n),p>n/m$, then $\psi$ has the same singularity type as $\phi$, more precisely:
$$\phi-C\Big(\|f\|_{L^{p}},p,\omega,\int_{X}\omega_{\phi}^{n}\Big)\leq\psi\leq\phi.$$
\end{thm}
\begin{proof}
Set
$$g(t):=(\mbox{Cap}_{m,\phi}\{\psi<\phi-t\})^{1/m}, t\geq0.$$
We will show that $g(M)=0$ for some $M$ under control. By Remark \ref{Cap m characterize m polar} we will then have $\psi\geq\psi-M$ a.e. with respect to $\omega^n$, which implies $\psi\geq\phi-M$ on $X$.

Since $H_m(\psi)=f\omega^n$, it follows from Proposition \ref{vol control cap} and Lemma \ref{1000}(take $v=\phi$) that
$$g(t+\delta)\leq\frac{C^{1/m}}{\delta}g(t)^{1+\alpha}, t>0, 0<\delta<1, \alpha>0.$$
Consequently, we can adapt \cite[Lemma 2.4]{EGZ09} to conclude that $g(M)=0$ for $M:=s_0+\frac{1}{1-2^{-\alpha}}$. As an important point, the constant $t_0>0$ has to be chosen so that
$$g(t_0)^{\alpha}<\frac{1}{2C^{1/m}}.$$

On the other hand, by Proposition \ref{107} and \cite[Corollary 3.19]{Lu13} we have $$g(t)\leq C_1\mbox{Cap}_{m}^{1/{m^2}}(\{\psi<-t\})\leq C_1/{t^{1/{m^2}}}.$$ 
We can take $t_0=C_1^{m^2/\alpha}2^{m^2/\alpha}C^{m/{\alpha}}$ and finish the proof.
\end{proof}
The next proposition tells us that if $\varphi,\psi$ are close in $m$-capacity, then they are close in $L^{\infty}$-norm, which has importance in the proof of our main theorem. It also generalizes \cite[Proposition 2.6]{EGZ09}.
\begin{prop}
Let $\varphi,\psi\in\mathcal{E}_{\phi}$ be two functions such that $\sup_X \varphi=\sup_X \psi=0$ and fix $\epsilon>0$. Assume that $H_m(\psi)=f\omega^n$ with $f\in L^p(X),p>n/m$ and $\varphi$ has the same singularity type as $\phi$. There exists  a constant $C>0$ such that
$$\sup_X (\varphi-\psi)\leq\epsilon+C[\mbox{Cap}_{m,\phi}(\psi-\varphi<-\epsilon)]^{\alpha/m}.$$

\end{prop}
\begin{proof}
 Just set $M:=\|\phi-\varphi\|_{L^{\infty}}$ and observe that $\varphi-\psi=\varphi-\phi-(\psi-\phi)$, meanwhile, Theorem \ref{relative boundedness} implies $\sup_X (\phi-\psi)<C\Big(\|f\|_{L^{p}},p,\omega,\int_{X}\omega_{\phi}^{n}\Big)$. Follow the line of \cite[Proposition 2.6]{EGZ09}.
\end{proof}
The following stability theorem was established in \cite{Lu13} for $\phi=0$.
\begin{prop}\label{1001}
Assume $H_m(\varphi)=f\omega^n,H_m(\psi)=g\omega^n$, where $\varphi,\psi\in\mathcal{E}_{\phi}$ and $f,g\in L^p(X)$ with $p>n/m$. Fix $r>0$. Then if $\gamma$ is taken so small that $\frac{\gamma mq}{r-\gamma(r+mq)}<\frac{mp-n}{(n-m)p},$ we have
$$\|\varphi-\psi\|_{L^{\infty}(X)}\leqslant C\|\varphi-\psi\|_{L^r(X)}^{\gamma},$$
where $q=\frac{p}{p-1}$ denotes the conjugate exponent of $p$, and the constant $C$ depends only on $n,m,p,r$ and upper bounds of $\|f\|_p,\|g\|_p$.
\end{prop}

Now we can prove our main result. We will give an alternative proof of special case of \cite{AAG21}. But as a byproduct, we obtain the relative boundedness of the solution of Hessian type equations, which is not contained in \cite{AAG22}.
\begin{lem}
Assume that $\phi\in\mbox{SH}_m(X,\omega)$. There exists a constant $C>0$ such that for all $\varphi\in\mbox{SH}_m(X,\omega)$ satisfying $\sup_X \varphi=0$, we have
$$\int_X (\varphi-\phi)\omega^n\geq -C.$$
It then follows that
$$\mathcal{C}:=\left\{\varphi\in SH_{m}(X,\omega)\big |\sup_{X}\varphi\leqslant0;\int_{X}(\varphi-\phi)\omega^{n}\geqslant-C_{0}\right\} $$
is a convex compact subset of $L^1(X)$.
\end{lem}

\begin{thm}
Assume that $\phi$ is a model potential. Let $F:X\times\mathbf{R}\rightarrow[0,+\infty)$ be a function satisfying the following conditions:
\begin{itemize}
\item[(i)] for almost every $x\in X, t\to F(x,t)$ is non-decreasing and continuous,
\item[(ii)] for any fixed $t\in\mathbf{R}$ there exists $p>n/m$ such that the function $x\to F(x,t)$ belongs to $L^p(X)$,
\item[(iii)] there exists $t_0\in\mathbf{R}$ such that $\int_X F(\cdot,t_0)=\int_X \omega^n$.
\end{itemize}
Then there exists a  unique  function $\varphi\in\mbox{SH}_m(X,\omega)$ up to a constant satisfying $[\phi]=[\varphi]$ and
$$\left(\omega+dd^{c}\varphi\right)^{m}\wedge\omega^{n-m}=F(x,\varphi)\omega^{n}.$$

\end{thm}
\begin{proof}
We adapt the idea of \cite{Lu13} with some necessary modification. We just give a sketch here.
To prove the existence and uniqueness, we take into two steps.

\noindent{\bf Case 1}. There exists $t_1\in\mathbb{R}$ such that $\int_{X}F(x,t_{1})\omega^{n}>\int_{X}F(x,t_{0})\omega^{n}.$

Take $\psi\in\mathcal{C}$, there exists $\varphi\in\mathcal{E}_{\phi}$ such that $\varphi$ has the same singularity type as $\phi$ and
\begin{equation}\label{eq proof}
H_m(\varphi)=F(\cdot,\psi+c_{\psi})\omega^n, \sup_X \varphi=0
\end{equation}
where $c_{\psi}\geq t_0$ is a constant such that
$$\int\limits_{X}F(.,\psi+c_{\psi})\omega^{n}=\int\limits_{X}\omega^{n}.$$ 
By a similar argument we can well define the map $\Phi:\mathcal{C}\rightarrow\mathcal{C},\psi\mapsto\varphi.$

Then we need to prove that $\Phi$ is continuous on $\mathcal{C}$. Assume that $(\psi_j)$ is a sequence in $\mathcal{C}$ converging to $\psi\in\mathcal{C}$ in $L^1(X)$ and let $\varphi_j=\Phi(\psi_j)$. Let $c_j:=c_{\psi_j}$ and we can prove that $(c_j)$ is uniformly bounded from the almost same argument. As a consequence, the sequence $(F(\cdot,\psi_j+c_j))_j\leq(F(\cdot,c_j))_j$ is uniformly bounded in $L^p(X),p>n/m$. By Theorem \ref{relative boundedness}, $\varphi_j-\phi$ is uniformly bounded.

Now we need to show that every cluster point of $(\varphi_j)$ satisfies $\Phi(\psi)=\varphi$ . Suppose that $\varphi_j\to\varphi$ in $L^1(X)$. Since $\varphi_j-\phi$ is uniformly bounded, it follows from Proposition \ref{1001} that $-C<\varphi-\phi<0$ for a constant $C$. By subtracting a subsequence if necessary we can assume that $\psi_j\to\psi$ almost everywhere on $X$ and $c_j\to c$. Since $t\mapsto F(x,t)$ is continuous we see that $F(\cdot,\psi_j+c_j)\to F(\cdot,\psi+c)$ almost everywhere. Thus $H_m(\varphi)=F(\cdot,\psi+C)$ which implies $\Phi(\psi)=\varphi$ and hence $\Phi$ is continous on $\mathcal{C}$.

Using the Schauder fixed point theorem, $\Phi$ has a fixed point in $\mathcal{C}$, denoted by $\varphi$. And we have $$H_m(\varphi)=F(\cdot,\varphi+c_{\varphi})\omega^n,$$
such that $\varphi$ has the same singularity type as $\phi$. The function $\varphi+c_{\varphi}$ is the  solution that we want.

\noindent{\bf Case 2.} $\int_X F(\cdot,t)\omega^n=\int_X F(\cdot,t_0)\omega^n, \forall t\geq t_0$.

For the second step, we just need to replace $C^{\prime}$ in \cite{Lu13} there by
$$\mathcal{C}^{\prime}:=\{\varphi\in\mbox{SH}_m(X,\omega)|-C_1\leq\varphi-\phi\leq0\}.$$

Similarly take $\psi\in\mathcal{C}^{\prime}$, we can find $\varphi\in\mathcal{E}_{\phi}$ such that $\sup_X \varphi=0$ and $$H_m(\varphi)=F(\cdot,\psi+c_{\psi})\omega^n,$$ where $t_0\leq c_{\psi}\leq t_0+C_1$ is constant such that $$\int_X F(\cdot,\psi+c_{\psi})\omega^n=\int_X \omega^n.$$
Then we can well define a continuous map $\Phi:\mathcal{C}\to\mathcal{C}^{\prime}$ by setting $\Phi(\psi)=\varphi$.

As in case 1, we can also  assume that $\psi_j\to\psi$ in $L^1(X)$. By Proposition \ref{1001} that the sequence $(\varphi_j)$ converges to $\varphi$ uniformly and $\varphi\in\mbox{SH}_m(X,\omega,\phi)$. By substracting a subsequence we can assume that $\psi_j\to\psi$ in $L^1(X)$ and $c_j\to c$. Then  we have that $H_m(\varphi)=F(\cdot,\psi+c)\omega^n$ and $\Phi(\psi)=\varphi$. After that, the continuity follows.

Applying the Schauder fixed point theorem, we have that $\Phi$ has a fixed point in $\mathcal{C}^{\prime}$, denoted by $\varphi\in\mathcal{E}_{\phi}$, having the same singularity type as $\phi$ and
$$H_m(\varphi)=F(\cdot,\varphi+c_{\varphi})\omega^n.$$
Finally, the function $\varphi+c_{\varphi}$ is the solution.
\end{proof}
\begin{rem}
For $\phi=0$ this theorem is proved in \cite{Lu13}, where the solution is continuous on $X$. But if $\phi$ has more singularities (for example $\phi$ is a model potential), we can not expect that the solution is continuous even bounded, as shown above. There are plenty of model potentials such as $\frac{1}{2}\omega$-psh functions with analytic singularities (see \cite[Proposition 5.23]{DNL23}\cite[Proposition 4.36]{DNL18}). This fact essentially used the resolution of Demailly's strong openness conjecture \cite{Dem} due to \cite{GZh15}.

As also shown in \cite{Lu13}, if $\phi=0$ and moreover $t\mapsto F(x,t)$ is strictly increasing for every $x\in X$, the solution is unique by using the continuity of solution.

\end{rem}

\section{Finite energy range of the Hessian operator}
In this chapter we will characterize the Borel measures $\mu$ that are equal to the Hessian product of some $u\in\mathcal{E}_{\chi}(X,\omega,m,\phi)$, where $\chi$ satisfies some natural polynomial growing condition and $\phi$ is a model potential $(\phi=P[\phi],\int_X H_m(\phi)>0)$.

$\chi$ is called a weight if it is a continuous increasing function from $[0,+\infty)$ to $[0,\infty)$ such that $\chi(0)=0,\chi(+\infty)=+\infty$ and satisfies the following condition
\begin{equation}
\forall t\geq 0,\forall \lambda\geq 1, \chi(\lambda t)\leq\lambda^M\chi(t),
\end{equation}
where $M$ is a fixed constant satisfying $1\leq M<\frac{n}{n-m}$.

Fix $\phi$ a model potential and set $\mathcal{E}_{\chi}(X,\omega,m,\phi)$ as the set of all $u\in\mathcal{E}(X,\omega,m,\phi)$ such that $E_{\chi}(u,\phi):=\int_X \chi(|u-\phi|)H_m(u)<+\infty$. For simplicity of notation, we denote $\mathcal{E}(X,\omega,m)=\mathcal{E}(X,\omega,m,0)$,
$\mathcal{E}_{\chi}(X,\omega,m)=\mathcal{E}_{\chi}(X,\omega,m,0)$ and $E_{\chi}(u)=E_{\chi}(u,0)$.

This section is a generalization of  \cite[\S6]{DNL23}. The proof is similar to the Monge-Amp\`ere case except for \cite[Proposition 6.11]{DNL23} because of the lack of integrability  theorem of $\omega-m$-subharmonic function. As a simple observation, we instead use the method of capacity estimate to overcome the difficulty there. This is the largest difference between our proof with Darvas-Nezza-Lu's.

For the reader's convenience, we follow the line of \cite[\S 6]{DNL23}.
\begin{lem}\label{const keep E}
If $u\in\mathcal{E}_{\chi}(X,\omega,m,\phi)$, $u+C\in\mathcal{E}_{\chi}(X,\omega,m,\phi)$ for any constant $C$.
\end{lem}
\begin{proof}
Since $\chi$ is increasing and being a weight, we have
\begin{align*}
\int_X \chi(|u+C-\phi|)H_m(u)&\leq\int_X \chi(|u-\phi|+|C|)H_m(u)\\
&=\int_X \chi(2(|u-\phi|+|C|)/2)H_m(u)\\
&\leq\int_X 2^M \chi((|u-\phi|+|C|)/2)H_m(u)\\
&\leq\int_X 2^M \chi(\max(\chi(|u-\phi|)),\chi(|C|))H_m(u)\\
&\leq 2^M\max(\int_X \chi(|u-\phi|),\chi(|C|) H_m(u))<+\infty
\end{align*}

See the Monge-Amp\`ere case \cite[Lemma 6.1]{DNL23}.
\end{proof}
\begin{lem}
There exist a uniform constant $C>0$ such that
$$\int_X \chi(\phi-u)H_m(\phi)\leq C$$ where $u$ belongs to $\mbox{PSH}(X,\omega,m,\phi)$ normalized by $\sup_X u=0$.
\end{lem}
\begin{proof}
The proof is similar to \cite[Lemma 6.2]{DNL23}.
By definition of weight $\chi$, $\chi(\phi-u)\leq(\phi-u)^M\chi(1)$ if $\phi-u\geq1$; $\chi(\phi-u)\leq\chi(1)$ if $\phi-u<1$. Hence
\begin{align*}
\int_X \chi(\phi-u)H_m(\phi)&\leq C^{\prime}\int_X (|\phi|^M+|u|^M+1)H_m(\phi)\\
&\leq C^{\prime}\int_X (|\phi|^M+|u|^M+1)\omega^n\leq C
\end{align*}
where the second inequality follows by $H_m(\phi)\leq\omega^n$ (by definition of Hessian product) and the last inequality follows from the fact that $\int_X |h|^M \omega^n$ is uniformly bounded for $h\in\mbox{PSH}(X,\omega,m)$ with $\sup_X h=0$. Indeed, 
\begin{align*}
\int_X |h|^M \omega^n &=\int_1^{+\infty}Mt^{M-1}\omega^n(h<-t)dt+\int_0^1 Mt^{M-1}\omega^n(h<-t)dt\\
&\leq C\int_1^{+\infty} Mt^{M-1}\mbox{Cap}_m^{\beta}(h<-t)+MVol(X)\\
&\leq C\int_1^{+\infty} Mt^{M-1-\beta}+MVol(X)\leq C
\end{align*}
if we choose $\beta\in (1,\frac{n}{n-m})$ such that $\beta=M+\epsilon$. Note that the first inequality follows from Lemma \ref{V geq Cap} for $\phi=0$ there. Note also that the second inequality follows from the fact that  there exists a constant C such that $\mbox{Cap}_m(f<-t)\leq \frac{C}{t}$ for all $\omega-m$-subharmonic function $f$ satisfying $\sup_X f=-1$ (see \cite[Corollary 3.19]{Lu13}).

\end{proof}
\begin{lem}
Let $u\in\mathcal{E}(X,\omega,m,\phi)$ satisfying $\sup_X u=0$, we have 
$$\int_X \chi(\phi-u)(\omega+dd^c u)^j\wedge(\omega+dd^c \phi)^{m-j}\wedge\omega^{n-m}\leq\int_X \chi(\phi-u)H_m(u).$$

\end{lem}
\begin{proof}
First note that we have $u\leq\phi\leq0$. Then use Corollary \ref{compare principle} and we will obtain the inequality, see \cite[Lemma 6.3]{DNL23}
\end{proof}
\begin{lem}
If $u,v\in\mathcal{E}(X,\omega,m,\phi)$ and $u,v\leq0$, we have
$$\int_X \chi(\phi-u)H_m(v)\leq 2^{m+M}E_{\chi}(u,\phi)+E_{\chi}(v,\phi).$$
\end{lem}
\begin{proof}
Note that $u\leq\phi\leq0$. The proof is essentially achieved by Corollary \ref{compare principle} and the observation that
$$\{2u\geq v+\phi-\chi^{-1}(t)\}\cap\{v\geq\phi-\chi^{-1}(t)\}\subseteq\{u\geq\phi-\chi^{-1}(t)\}.$$ 
See \cite[Lemma 6.4]{DNL23}.

\end{proof}
\begin{lem}\label{ineq of energy}
If $u,v\in\mathcal{E}(X,\omega,m,\phi)$ satisfy $u\leq v\leq0$, then
$$\int_{X}\chi(\phi-v)H_m(v)\leq\int_{X}\chi(\phi-u)H_m(v)\leq2^{n+M}E_{\chi}(u,\phi).$$
\end{lem}
\begin{proof}
Note that we have $u\leq v\leq\phi$ for the same reason above. Use this fact and play the same trick we will obtain the result. See \cite[Lemma 6.5]{DNL23}.
\end{proof}
\begin{prop}
If $u\in\mathcal{E}_{\chi}(X,\omega,m,\phi)$ and $u\leq v$, we have $v\in\mathcal{E}_{\chi}(X,\omega,m,\phi)$. Moreover, $\mathcal{E}_{\chi}(X,\omega,m,\phi)$ is convex.
\end{prop}
\begin{proof}
The proof is similar to \cite[Proposition 6.6]{DNL23}. We may assume that $v\leq0$ by Lemma \ref{const keep E}. The Lemma above then implies the first statement. We only need to prove the second statement. If $u,v\in\mathcal{E}_{\chi}(X,\omega,m,\phi)$, then by \cite[Proposition 3.22]{LN22} $P(u,v)\in\mathcal{E}(X,\omega,m,\phi)$. From \cite[Corollary 3.11]{LN22} we have
\begin{align*}
\int_{X}\chi(\phi-P(u,v))H_m(P(u,v))&\leq\int_{\{P(u,v)=u\}}\chi(\phi-u)H_m(u)+\int_{\{P(u,v)=v\}}\chi(\phi-v)H_m(v)\\
&\leq\int_{X}\chi(\phi-u)\theta_{u}^{n}+\int_{X}\chi(\phi-v)\theta_{v}^{n}<\infty.
\end{align*}
Therefore $P(u,v)\in\mathcal{E}_{\chi}(X,\omega,m,\phi)$. But as a obvious observation that $tu+(1-t)v\geq w$, the result holds by Lemma \ref{ineq of energy}.
\end{proof}
\begin{lem}\label{L1 convergence energy}
Assume $(u_j)$ is a sequence in $\mathcal{E}_{\chi}(X,\omega,m,\phi)$ converging in $L^1$ to $u\in\mbox{PSH}(X,\omega,\phi)$. If $\sup_j E_{\chi}(u_j,\phi)<+\infty$, then $u\in\mathcal{E}_{\chi}(X,\omega,m,\phi)$.
\end{lem}
\begin{proof}
The proof is similar to \cite[Lemma 6.7]{DNL23}, so we only give a sketch. By Lemma \ref{ineq of energy} we can assume that $u_j$ decrease to $u$. Fix $t>0$ and set $u_{j,t}:=\max(u_j,\phi-t)$, then $u_{j,t}$ decrease to $u_t:=\max(u,\phi-t)$ as $k\to\infty$. By Lemma \ref{ineq of energy} $E_{\chi}(u_{j,t},\phi)$ is uniformly bounded with respect to $t$. Note that $\chi(\phi-u_{j,t})$ are quasi-continuous and uniformly bounded. Applying \cite[Theorem 3.3]{LN22} we have
$$
\chi(t)\int_{\{u\leq\phi-t\}}H_m(u_t)\leq\int_X \chi(\phi-u_t)H_m(u_{j,t})\leq\liminf_{k\to\infty}\int_X \chi(\phi-u_{j,t})H_m(u_{j,t})\leq C.$$
Since $\chi(t)\to\infty$ as $t\to\infty$, the integral $\int_{\{u\leq\phi-t\}}H_m(u_t)\to0$. Using the plurifine property of the Hessian operator we will get that $u\in\mathcal{E}(X,\omega,m,\phi)$. Note also that
$$\int_{\{u>\phi-t\}}\chi(\phi-u_t)H_m(u_t)=\int_{\{u>\phi-t\}}\chi(\phi-u)H_m(u)\leq C.$$
Let $t\to\infty$ we have $u\in\mathcal{E}_{\chi}(X,\omega,m,\phi)$.

\end{proof}
\begin{lem}
Let $\mu$ be a positive Borel measure on $X$. Assume that $\mu\{\phi=-\infty\}=0$ and $\chi(\phi-u)\in L^1,\forall u\in\mathcal{E}_{\chi}(X,\omega,m,\phi)$. Fix a constant $A>0$. Then there exists a constant $C>0$ depending on $A$ such that for all $u\in\mathcal{E}(X,\omega,m,\phi)$ satisfying $\sup_X u=0$ and $E_{\chi}(u,\phi)\leq A$ we have
$$\int_X \chi(\phi-u)d\mu\leq C.$$
\end{lem}
\begin{proof}
We sketch the similar proof as \cite[Lemma 6.8]{DNL23}. Assume by contradiction that there exists a sequence $(u_j)\in\mathcal{E}_{\chi}(X,\omega,m,\phi)$ such that $\sup_X u_j=0$ and $E_{\chi}(u_j,\phi)\leq A$ such that
$$\int_X \chi(\phi-u_j)d\mu\geq 4^{jM}.$$ 
Set $v_k:=P(\min_{1\leq j\leq k}(2^{-j}u_j+(1-2^{-j}\phi)))\leq\phi$. We can prove that $\int_X \chi(\phi-v_k)H_m(v_k)$ has a uniform upper bound. Consider the decreasing limit of $v_k$ denoted by $v$. By Lemma \ref{L1 convergence energy} we have $v\in\mathcal{E}_{\chi}(X,\omega,m,\phi)$. But $\int_X \chi(\phi-u)d\mu=\infty$, which makes a contradiction.
\end{proof}
\begin{lem}\label{int leq weighted energy}
 Assume that $\mu\{\phi=-\infty\}=0$ and $\chi(\phi-u)\in L^1,\forall u\in\mathcal{E}_{\chi}(X,\omega,m,\phi)$.
Then there exists a constant $C>0$ such that for all $u\in\mathcal{E}(X,\omega,m,\phi)$ satisfying $\sup_X u=0$, we have
$$\int_X \chi(\phi-u)d\mu\leq C(E_{\chi}(u,\phi)+1).$$ 
\end{lem}
\begin{proof}
The proof is similar to \cite[Lemma 6.9]{DNL23}.
\end{proof}
\begin{lem}\label{Lemma 6.10}
There exists a constant $C>0$ such that $\forall u,v\in\mathcal{E}_{\chi}(X,\omega,m,\phi)$ with $\sup_X v=0$ and $u\leq0$, we have
$$\int_X \chi(\phi-v)H_m(u)\leq C(1+E_{\chi}(u,\phi))E_{\chi}(v,\phi)^{M/(M+1)}+C.$$
\end{lem}
\begin{proof}
The proof is similar to \cite[Lemma 6.10]{DNL23}\cite{DV21}.
\end{proof}
\begin{lem}\label{int leq energy}
If $\mu\leq A\mbox{Cap}_{m,\phi}$ for some constant $A>1$. Then there exists a constant $B>0$ depending on $A$ such that for all $u\in\mbox{SH}_m(X,\omega,m,\phi)$ with $\sup_X u=-1$ we have
$$\int_X (\phi-u)^2d\mu\leq B(E(u,\phi)+1)^r,$$
where $r$ is a large constant.
\end{lem}
\begin{proof}
The proof is motivated by \cite[Lemma 4.18]{DNL18}(essentially goes back to \cite[Lemma 2.9]{BBGZ13}) and \cite[Lemma 6.8]{LN15}. We can assume that $u\in\mathcal{E}_{\chi}(X,\omega,m,\phi)$. For each $t>1$ we set $u_t:=t^{-1}u+(1-t^{-1})\phi$. Fix $\psi\in\mbox{SH}_m(X,\omega,m)$ satisfying $-1\leq\psi-\phi\leq0$. Then we know $u_t,\psi\in\mathcal{E}(X,\omega,m,\phi)$ and the following inclusion 
$$(u<\phi-2t)\subset(u_{t}<\psi-1)\subset(u<\phi-t),t>1.$$
By comparison principle we have
$$H_m(\psi)(u<\phi-2t)\leq H_m(\psi)(u_t<\psi-1)\leq H_m(u_t)(u_t<\psi-1)\leq H_m(u_t)(u<\phi-t).$$
Now we have
\begin{align*}
\int_{1}^{\infty}tH_m(\psi)(u<\phi-2t)dt&\leq\int_{1}^{\infty}tH_m(u_t)(u<\phi-t)dt\\
&\leq\int_{1}^{\infty}tH_m(\phi)(u<\phi-t)dt+\sum_{k=1}^{m}{\binom{m}{k}}\int_1^{\infty}\omega_{u}^{k}\wedge\omega_{\phi}^{m-k}\wedge\omega^{n-m}(u<\phi-t)\\
&\leq\int_{1}^{\infty}tH_m(\phi)(u<\phi-t)dt+CE(u,\phi).
\end{align*}

The last inequality above used comparison principle. We now consider the first term in the last inequality.
\begin{align*}
\int_{1}^{\infty}tH_m(\phi)(u<\phi-t)dt&\leq\int_1^{\infty}t Vol(u<\phi-t)=\int_1^{\infty}Vol(u<\phi-t)^{\gamma}Vol(u<\phi-t)^{1-\gamma}dt\\
&\leq\left[\int_{1}^{+\infty}t\operatorname{Vol}(u<\phi-t)^{q\gamma}\mathrm{d}t\right]^{1/q}\left[\int_{1}^{+\infty}t\operatorname{Vol}(u<-t)^{r(1-\gamma)}\mathrm{d}t\right]^{1/r}\\
&\leq A^{\prime\prime}\bigg[\int_{1}^{+\infty}t\mathrm{Cap}_{m,\phi}(u<\phi-t)^{pq\gamma}\mathrm{d}t\bigg]^{1/q}\times\left[\int_{1}^{+\infty}t\mathrm{Cap}_{m}(u<-t)^{pr(1-y)}\mathrm{d}t\right]^{1/r}\\
&\leq A^{\prime}\bigg[\int_{1}^{+\infty}t\mathrm{Cap}_{m,\phi}(u<\phi-2t)\mathrm{d}t\bigg]^{1/q}\bigg[\int_{1}^{+\infty}t^{1-pr(1-p)}\mathrm{d}t\bigg]^{1/r}.
\end{align*}
Note that the constant $A^{\prime}$ does not depend on the choice of $u$ because of the normalization condition. Here $\frac{1}{q}+\frac{1}{r}=1$ and $p$ can be arbitrarily chosen in $(1,n/(n-m))$. We can  also choose $\gamma$ so that $pq\gamma=1$ and $pr(1-\gamma)>2$. Note that the first inequality follows from \cite[Corollary 3.11]{LN22} and the third inequality follows from Lemma \ref{V geq Cap}. Note also that the last inequality follows from \cite[Corollary 3.19]{Lu13}.

Set $u_j:=\max(u,\phi-j)$ and replace $u$ by $u_j$ then we will get
$$C_j\leq A^{\prime}C_j^{1/q}+C\cdot E(u_j,\phi),$$ where $A^{\prime},C>0$ is a constant and $C_j:=\int_1^{\infty}t\mbox{Cap}_{m,\phi}(u_j<\phi-2t)$. It implies that
$C_j\leq C(1+E(u,\phi))^r$ for some $C>1$ and $r>1$. Finally, we can wirte
\begin{align*}
\int_X (\phi-u_j)^2d\mu&=2\int_0^{\infty}t\mu(u_j<\phi-t)dt\leq \mu(X)+2\int_1^{\infty}t\mu(u_j<\phi-t)dt\\
&\leq \mu(X)+8\int_1^{\infty}At\mbox{Cap}_{m,\phi}(u_j<\phi-2t)dt\leq B(E(u_j,\phi)+1)^r.
\end{align*}
By the monotone convergence theorem and the property of Hessian operator, the result holds.
\end{proof}
\begin{rem}
The reader can compare the above lemma with Lemma \ref{int leq weighted energy}. By Lemma \ref{int leq energy} we know $(\phi-u)^2\in L^1(\mu),\forall u\in\mathcal{E}_{\chi}$ if $\mu$ is dominated by $\mbox{Cap}_{m,\phi}$. By the proof of Lemma \ref{int leq weighted energy}, there exists a constant $C>0$ such that $\int_X (\phi-u)^2d\mu\leq C(E(u,\phi)+1),\forall u\in\mbox{SH}_m(X,\omega,\phi),\sup_X u=0$.

\end{rem}

\begin{thm}\label{M_A solution}
Assume that $\mu\leq A\mbox{Cap}_{m,\phi}$ for some $A\geq1$. Then there exists uniquely $u\in\mathcal{E}^1(X,\omega,m,\phi),\sup_X u=-1$ such that $\mu=H_m(u)$.
\end{thm}
\begin{proof}
The proof is similar to \cite[Theorem 4.25]{DNL18}.
\end{proof}
\begin{cor}\label{special Ep solution}
Let $\mu\leq A \mbox{Cap}_{m,\phi}$ for some $A\geq1$. Then there exists a unique $\omega$-$m$-subharmonic function $\psi\in\bigcap_{p\geq1}\mathcal{E}_{m,\phi}^{p},\sup_X \psi=-1$ such that $\mu=H_m(\psi)$.
\end{cor}
\begin{proof}
By the proof of Lemma \ref{int leq energy} we have that for all $u\in\mbox{SH}_m(X,\omega,\phi)$ with $\sup_X u=-1$, $$\int_X (\phi-u)^pd\mu\leq B_p(E^{p-1}(u,\phi)+1)^r,$$
where $B_p$ only depends on $p$. Now by Theorem \ref{M_A solution} there exists a unique $\omega$-$m$-subharmonic function $\psi\in\mathcal{E}^1(X,\omega,m,\phi),\sup_X \psi=-1$ such that $\mu=H_m(\psi)$. Applying to $u:=\psi$ and by induction we will obtain the result. See also the Monge-Ampere case \cite[Theorem 3.7]{DV21}.
\end{proof}
\begin{prop}\label{6.13}
Assume that $\mu$ is a positive measure satisfying $\int_X H_m(\phi)=\mu(X)>0$ and $\mu(\{\phi=-\infty\})=0$. Assume also that
$$\int_{X}\chi(\phi-\varphi)d\mu\leq aE_{\chi}(\varphi,\phi)+C,\quad\varphi\in\mathcal{E}_{\chi}(X,\omega,m,\phi),\sup_{X}\varphi=0,$$
for some constant $a\in (0,1),C>0$. Then we have $\mu=H_m(u)$ for some $u\in\mathcal{E}_{\chi}(X,\omega,m,\phi)$.
\end{prop}
\begin{proof}
The proof is inspired by \cite[Proposition 6.13]{DNL23}. First we claime that $\mu$ does not charge $m$-polar sets. If $E$ is a $m$-polar Borel subset, by the proof of \cite[Lemma 6.12]{DNL23} we have that $E\subset\{h=-\infty\}$ for some $h\in\mathcal{E}_{\chi}$. Since $\int_E \chi(\phi-u)d\mu\leq\int_X \chi(\phi-u)d\mu<\infty$, we get that $\mu$ does not charge $E\cap\{\phi\neq-\infty\})$. Then it follows that $\mu(E)=0$ because $\mu(\{\phi=-\infty\})$. By the similar arguement as \cite[Lemma 4.26]{DNL18}(originally goes back to \cite{Ce98}), there exists $\nu\leq A\mbox{Cap}_{m,\phi},A\geq1$ and $0\leq f\in L^1(X,\nu)$ such that $\mu=f\nu$. Now for each $j>1$ we let $\varphi_j$ be the unique solution of $$H_m(\varphi_j)=c_j\min(f,j)\nu,\sup_X \varphi_j=0,$$ where $c_j$ is a constant to have equality between the total masses of the left and right hand side and $\varphi_j\in\mathcal{E}^1(X,\omega,m,\phi)$. Then we have $E_{\chi}(\varphi_j,\phi)$ is finite. Indeed, by Corollary \ref{special Ep solution}(replace $\psi$ there by $\varphi_j-1$)
$$E_{\chi}(\varphi_j,\phi)\leq\int_X \chi(\phi-\varphi_j+1)H_m(\varphi_j)\leq\int_X (\phi-\varphi_j+1)^{[M]+1}\chi(1)\max(1,c_jj)d\nu<\infty.$$ We claim that this bounded in uniform in $j$.

Indeed, since $H_m(\varphi_j)\leq c_j fd\nu=c_jd\mu$, we have
$$E_{\chi}(\varphi_j,\phi)\leq\int_X \chi(\phi-\varphi_j)c_jd\mu\leq ac_jE_{\chi}(\varphi_j,\phi)+C$$
implies $E_{\chi}(\varphi_j,\phi)<C(1-\lambda)^{-1}$ where $c_ja<\lambda<1$. This is possible because $c_j\to1$ when $j\to\infty$. By extracting a subsequence if necessary we can assume that $\varphi_j\to\varphi$ in $L^1$. By Lemma \ref{L1 convergence energy} we have $\varphi\in\mathcal{E}_{\chi}$. By the proof of \cite[Lemma 5.16]{DNL23} we have $H_m(\varphi)\geq\mu$. Comparing the total mass, we get the equality and the result holds.
\end{proof}
We can now prove our main result.
\begin{thm}
Assume $\mu$ is a Radon measure such that $\mu(\{\phi=-\infty\})=0$ and $\int_X H_m(\phi)=\mu(X)>0$. Then the following statements are equivalent:
\begin{itemize}
\item[(i)] There exists a constant $C>0$ such that for all $u\in\mathcal{E}_{\chi}(X,\omega,m,\phi)$ with $\sup_X u=0$, we have
$$\int_{X}\chi(\phi-u)d\mu\leq CE_{\chi}(u,\phi)^{M/(M+1)}+C.$$
\item[(ii)]$\chi(|\phi-u|)\in L^{1}(\mu),\forall u\in\mathcal{E}_{\chi}(X,\omega,m,\phi)$.
\item[(iii)]$\mu=H_m(\varphi)$ for some $\varphi\in\mathcal{E}_{\chi}(X,\omega,m,\phi)$ satisfying $\sup_X \varphi=0$.

\end{itemize}
\end{thm}
\begin{proof}
$\mathrm{(i)}\implies\mathrm{(ii)}$ is obvious. Now we prove $\mathrm{(ii)}\implies\mathrm{(iii)}$. The proof is inspired by \cite[Theorem 6.14]{DNL23}. Let $v\in\mathcal{E}(X,\omega,m,\phi)$ be the unique solution to $H_m(v)=(4C)^{-1}\mu+b\omega^n,\sup_X v=-1$. Here $b>0$ is a constant such that $\int_X (4C)^{-1}\mu+b\omega^n=\int_X H_m(\phi) $ and $C$ is the constant in Lemma \ref{int leq weighted energy}. Then by Lemma \ref{int leq weighted energy} there exists a constant $C_1>0$ such that for all $\varphi\in\mathcal{E}_{\chi}(X,\omega,m,\phi)$ with $\sup_X \varphi=0$ we have
$$\int_X \chi(\phi-u)H_m(v)\leq (4C)^{-1}\int_X \chi(\phi-u)d\mu\leq \frac{1}{4}E_{\chi}(\varphi,\phi)+C_1.$$
By Proposition \ref{6.13} we have $v\in\mathcal{E}_{\chi}$. Applying Lemma \ref{Lemma 6.10} we have that 
$$
\int_X \chi(\phi-u)d\mu\leq 4C\int_X \chi(\phi-u)H_m(v)\leq C(1+E_{\chi}(v,\phi))E_{\chi}(u,\phi)^{M/(M+1)}+C.$$
We can take a large number $N>0$ such that $(E_{\chi}(u-N,\phi))^{1/(M+1)}\geq 2^{M+1}C(1+E_{\chi}(v,\phi))$ in above inequality. This is possible because of the definition $\chi$. As a result,
$$\int_X \chi(\phi-u)d\mu\leq\int_X \chi(\phi-u+N)d\mu\leq \frac{1}{2}E_{\chi}(u,\phi)+C,$$
where the last inequality follows from the same trick of Lemma \ref{const keep E}.
Applying Proposition \ref{6.13} we can finish the argument. At last, Lemma \ref{Lemma 6.10} tells us $\mathrm{(iii)}\implies\mathrm{(i)}$. 
\end{proof}

\noindent{\bf Acknowledgement}

I would like to deeply thank my advisor, Prof. Xiangyu Zhou for his consistent guidance for me and providing me interesting problems at his seminar.


\begin{thebibliography}{}

\bibitem[AAG22]{AAG22}
H. Amal,S. Asserda and A. El-Gasmi, Weak Solutions to Complex Hessian Type Equations in the Class $\mathcal{E}_{\phi}(X,\omega,m)$, Volume 52, pages 117–128, 2022.

\bibitem[AV10]{AV10}
S. Alesker and M. Verbitsky, Quaternionic Monge-Amp\`ere equation and Calabi problem for HKT-manifolds, Israel J. Math. 176 (2010), 109-138.

\bibitem[BBGZ13]{BBGZ13}
R. J. Berman, S. Boucksom, V. Guedj, and A. Zeriahi. A variational approach to complex
Monge-Amp\`ere equations, Publ. Math. Inst. Hautes Etudes Sci. 117 (2013), pp. 179–245.

\bibitem[BEGZ10]{BEGZ10}
 S. Boucksom, P. Eyssidieux, V. Guedj, A. Zeriahi. Monge-Amp\`ere equations in big cohomology
classes, Acta Math. (2010), Volume 205, Issue 2, pp 199–262.

\bibitem[Blo05]{Blo05}
Z. B\l ocki, Weak solutions to the complex Hessian equation, Ann. Inst. Fourier (Grenoble) 55 (5) (2005) 1735–1756.

\bibitem[Ce98]{Ce98}
U. Cegrell. “Pluricomplex energy”. In: Acta Math. 180.2 (1998), pp. 187–217

\bibitem[Dem]{Dem}
 J.-P. Demailly. Multiplier ideal sheaves and analytic methods in algebraic geometry. School on
Vanishing Theorems and Effective Results in Algebraic Geometry (Trieste, 2000), 1–148, ICTP
Lect. Notes, 6, Abdus Salam Int. Cent. Theoret. Phys., Trieste, 2001

\bibitem[DK14]{DK14}
2] S. Dinew and S. Ko\l odziej, A priori estimates for complex Hessian equations, Anal. PDE 7 (2014), no. 1, 227-244.

\bibitem[DK17]{DK17}
S. Dinew and S. Kolodziej, Liouville and Calabi-Yau type theorems for complex Hessian equations, Amer. J. Math. 139 (2017), no. 2, 403–415.

\bibitem[DNL18]{DNL18}
T. Darvas, E. Di Nezza, and C. H. Lu. Monotonicity of nonpluripolar products and complex
Monge-Amp\`ere equations with prescribed singularity, Anal. PDE 11.8 (2018), pp. 2049–
2087.

\bibitem[DNL21]{DNL21}
 T. Darvas, E. Di Nezza, and C. H. Lu. Log-concavity of volume and complex Monge-Amp`ere
equations with prescribed singularity, Math. Ann. 379.1-2 (2021), pp. 95–132.


\bibitem[DNL18]{DNL18}
T. Darvas, E.Di Nezza, and C.H. Lu. On the singularity type of full mass currents in big cohomology classes,
Compositio mathematica 154(2) (2018): 380–409.

\bibitem[DNL23]{DNL23}
T. Darvas, E.Di Nezza, and C.H. Lu. Relative Pluripotential Theory on Compact K\"ahler Manifolds,
arXiv:2303.11584v1.

\bibitem[DV21]{DV21}
D.-T. Do and D.-V. Vu. “Complex Monge-Amp`ere equations with solutions in finite energy
classes”, Math. Research Letters (2021).

\bibitem[EGZ09]{EGZ09}
 P. Eyssidieux, V. Guedj, and A. Zeriahi. Singular K\"ahler-Einstein metrics, J. Amer. Math. Soc. 22
(2009), no. 3, 607–639.

\bibitem[GLZ19]{GLZ19}
V. Guedj, C. H. Lu, and A. Zeriahi. Plurisubharmonic envelopes and supersolutions, J. Differential Geom. 113.2 (2019), pp. 273–313.

\bibitem[GZ05]{GZ05}
 V. Guedj and A. Zeriahi. Intrinsic capacities on compact K\"ahler manifolds, J. Geom. Anal.
15.4 (2005), pp. 607–639.

\bibitem[GZ07]{GZ07}
 V. Guedj and A. Zeriahi, The weighted Monge-Amp\`ere energy of quasi plurisubharmonic functions,
Journal of Functional Analysis, 250(2), 442–482.

\bibitem[GZ17]{GZ17}
V. Guedj and A. Zeriahi. Degenerate complex Monge-Amp\`ere equations. Vol. 26. EMS Tracts in
Mathematics. European Mathematical Society (EMS), Z\"urich, 2017, pp. xxiv+472.


\bibitem[GZh15]{GZh15}
Q. Guan, X. Zhou. A proof of Demailly’s strong openness conjecture, Annals of Math. (2015), Issue 2, Volume 182, Pages 605–616

\bibitem[HMW10]{HMW10}
 Z. HOU, X. N. Ma, and D.M. Wu, A second order estimate for complex Hessian equations on a compact K\"ahler manifold, Math. Res. Lett. 17 (2010), no. 3, 547-561.

\bibitem[Hou09]{Hou09}
 Z. Hou, Complex Hessian equation on Kähler manifold, Int. Math. Res. Not. IMRN (16) (2009) 3098–3111.

\bibitem[Jbi10]{Jbi10}
 A. Jbilou, Equations hessiennes complexes sur des variétés kählériennes compactes, C. R. Math. Acad. Sci. Paris 348 (1–2) (2010) 41–46.

\bibitem[Kol98]{Kol98}
 S. Ko\l odziej. The complex Monge-Amp\`ere equation, Acta Math. 180 (1998), no. 1, 69–117.

\bibitem[Kol03]{Kol03}
S. Ko\l odziej. The Monge-Amp\`ere equation on compact K\"ahler manifolds, Indiana Univ. Math. J. 52 (2003), 667–686.

\bibitem[Kok10]{Kok10}
 V.N. Kokarev, Mixed volume forms and a complex equation of Monge–Ampère type on K\"ahler manifolds of positive curvature, Izv. Ross. Akad. Nauk Ser. Mat. 74 (3) (2010) 65–78.

\bibitem[LN15]{LN15}
H.C. Lu and V.-D. Nguyen. Degenerate comlex Hessian equations on compact K\"ahler manifolds, Indiana University Mathematics Journal, 2015, Vol. 64, No. 6 (2015), pp. 1721-
1745.

\bibitem[LN22]{LN22}
H.C. Lu and V.-D. Nguyen. Complex Hessian equations with prescribed singularity on compact K\"ahler manifolds, with V.D. Nguyen, Ann. Sc. Norm. Super. Pisa Cl. Sci. (5) 23 (2022), no. 1, 425--462.

\bibitem[Lu13]{Lu13}
H.C. Lu. Solutions to degenerate complex Hessian equations , J. Math. Pures Appl. (9) 100 (2013), no. 6, 785--805.

\bibitem[Lu21]{Lu21}
H.C. Lu. Comparison of Monge-Amp\`ere capacities, Ann. Polon. Math. 126.1 (2021).

\bibitem[Ngu16]{Ngu16}
V. D. Nguyen. The complex Monge–Ampère type equation on compact Hermitian manifolds and applications,  Advances in Mathematics 286 (2016) 240–285.

\bibitem[PPZ17]{PPZ17}
Duong H. Phong, Sebastien Picard, and Xiangwen Zhang, The Fu-Yau equation with negative slope parameter, Invent. Math. 209 (2017), no. 2, 541–576.

\bibitem[PPZ18]{PPZ18}
D. H. Phong, S. Picard, and X. Zhang, Fu-Yau Hessian equations, J. Differ. Geom. 118 (2021), no. 1, 147–187

\bibitem[PPZ19]{PPZ19}
Duong H. Phong, Sebastien Picard, and Xiangwen Zhang, On estimates for the Fu-Yau generalization of a Strominger system, J. Reine Angew. Math. 751 (2019), 243–274

\bibitem[Pli13]{Pli13}
S. Plis. The smoothing of m-subharmonic functions, perprint, available at http://arxiv.org/abs/arXiv:l312.1906.

\bibitem[Wit19]{Wit19}
 D. Witt Nystr\"om. Monotonicity of non-pluripolar Monge-Amp\`ere masses, Indiana Univ. Math. J. 68.2 (2019), pp. 579–591.

\bibitem[Xia19]{Xia19}
M. Xia. Integration by parts formula for non-pluripolar product, In: arXiv:1907.06359 (2019).

\bibitem[Yau78]{Yau78}
 S. T. YAU, On the Ricci curvature of a compact Kähler manifold and the complex
 Monge-Ampère equation. I, Comm. Pure Appl. Math. 31 (1978), no. 3, 339-411.

\end{thebibliography}
		\end{document}